\documentclass{panm_mio}

\usepackage{graphicx}

\usepackage{amsthm}
\usepackage{amsmath,amsfonts}
 \usepackage{cite}

\def\Sm{\mathcal{C}_\Theta^{\infty}}
\def\D{\mathcal{D}(I)}

\title{The $\star$-product approach for linear ODEs: a numerical study of the scalar case}

\author{Stefano Pozza$^1$, Niel Van Buggenhout$^1$\\
\vskip 2mm {\small
$^1$
Charles University \\
Sokolovská 83 186, 75 Praha 8, Czech Republic \\
pozza@karlin.mff.cuni.cz, buggenhout@karlin.mff.cuni.cz \\
}
}

\abstract{
Solving systems of non-autonomous ordinary differential equations (ODE) is a crucial and often challenging problem. Recently a new approach was introduced based on a generalization of the Volterra composition. In this work, we explain the main ideas at the core of this approach in the simpler setting of a scalar ODE. Understanding the scalar case is fundamental since the method can be easily extended to the more challenging problem of systems of ODEs. Numerical examples illustrate the method's efficacy and its properties in the scalar case.
}



\begin{document}

\maketitle

\section{Introduction}
Systems of non-autonomous linear ordinary differential equations arise in a variety of contexts \cite{Autler1955,BenEtAll17,Blanes15,kwaSiv72,Lauder1986,Shirley1965}. Yet, their solution remains surprisingly difficult to obtain, both formally and numerically, especially when dealing with systems of large-to-huge size.
Consider an $N\times N$ matrix $\tilde{A}(t)$ depending on the variable $t\in I\subseteq \mathbb{R}$. The unique solution $U_s(t)$ of the system
\begin{equation}\label{eq:ode:intro} 
\tilde{A}(t) U_s(t)=\frac{d}{dt}U_s(t), \quad U_s(s)=I_N, \quad \text{ for } t \geq s, \quad t,s\in I,
\end{equation}
with $I_N$ the $N \times N$ identity matrix, is a $N \times N$ matrix-valued function known as the \emph{time-ordered exponential} of $\tilde{A}(t)$. If $\tilde{A}(\tau_1)\tilde{A}(\tau_2)=\tilde{A}(\tau_2)\tilde{A}(\tau_1)$ for all $\tau_1,\tau_2 \in I$, then the time-ordered exponential can be expressed as $$U_s(t)=\exp\left(\int_s^{t} \tilde{A}(\tau)\, \text{d}\tau\right).$$ 
In general, however, $U_s$ has no known simple expression in terms of $\tilde{A}(t)$.

In \cite{Giscard2015,BonGis20}, a new expression for the solution is given using the \emph{path-sum approach}, a method able to express each element of $U_s(t)$ as a finite sequence of integral equations. However, this requires solving an NP-hard problem. In \cite{ProceedingsPaper2020,GiscardPozza2021,GiscardPozza2022}, the NP-hard problem is overcome by introducing the \emph{$\star$-Lanczos method}, a constructive method able to tridiagonalize $\tilde{A}(t)$.
At the heart of both the path-sum and $\star$-Lanczos method is a non-commutative convolution-like product, denoted by $\star$, defined between certain distributions \cite{schwartz1978}. Thanks to this product, the solution of \eqref{eq:ode:intro} can be expressed through the $\star$-product inverse \cite{ProceedingsPaper2020}.

In this work, we aim to illustrate to the numerical mathematics community the $\star$-product and how it can be used to solve an ODE numerically. For this reason, we restrict the presentation to the simpler case in which the ODE \eqref{eq:ode:intro} is a scalar equation. While this framework may look too simple to show the potential of the newly introduced technique, the reader should keep in mind that the results and construction we illustrate for the scalar case can be straightforwardly extended to the matrix case in full generality.

In Section \ref{sec:starexpression}, we give an introduction to the $\star$-product and the related expression for the solution of a scalar ODE. Section \ref{sec:starsol} discretizes the $\star$-product. As a consequence, the ODE solution can be obtained by solving a linear system. Several properties of the linear system are numerically investigated in Section \ref{sec:matrices}. 
The numerical experiments in Section \ref{sec:numexp} show that the presented strategy can compute the solution up to machine precision.
Section \ref{sec:conclusion} concludes the presentation.

\section{ODE solution by the $\star$-product approach}\label{sec:starexpression}
Given two appropriate bivariate functions $\tilde{f}_1(t,s), \tilde{f}_2(t,s)$,
the \emph{Volterra composition}, introduced by Vito Volterra (e.g., \cite{Volterra1928}), is defined as
\begin{equation*}
  \big(\tilde{f}_2 \star_v \tilde{f}_1\big)(t,s) := \int_s^{t} \tilde{f}_2(t,\tau) \tilde{f}_1(\tau, s) \, \text{d}\tau.
\end{equation*}
For our purposes, it suffices to assume $\tilde{f}_1$ and $\tilde{f}_2$ to be smooth (i.e., infinitely differentiable) on both variables over a bounded interval $I = [0, T ]$  to have a well-defined operation for every $t, s \in I$. Therefore, from now on, a function marked with a tilde will stand for a smooth function in both $t$ and $s$ over $I$.
Since the Volterra composition is closed for such functions, we are allowed to define the $k$th $\star_v$-power of a function $\tilde{f}$, that is, $\tilde{f}^{\star_v 1} = \tilde{f}$, and 
\begin{align*}
    \tilde{f}^{\star_v k} &:= \tilde{f} \star_v \tilde{f} \dots \star_v \tilde{f} = \\
    &=  \int_s^{t} \tilde{f}(t,\tau_1) \int_s^{\tau_1} \tilde{f}(\tau_1,\tau_2) \, \dots \int_s^{\tau_{k-1}} \tilde{f}(\tau_{k-2},\tau_{k-1})\tilde{f}(\tau_{k-1},s)   \, \text{d}\tau_{k-1} \cdots \, \text{d}\tau_2  \; \text{d}\tau_1,    
\end{align*}
for $k \geq 1$, with the convention $\tau_0 = t$. 
Moreover, the operation is also defined for univariate functions $\tilde{f}_2(t)$:
\begin{equation*}
  \big(\tilde{f}_2(t) \star_v \tilde{f}_1(t,s)\big)(t,s) := \int_s^{t} \tilde{f}_2(t) \tilde{f}_1(\tau, s) \, \text{d}\tau =  \tilde{f}_2(t) \int_s^{t} \tilde{f}_1(\tau, s) \, \text{d}\tau.
\end{equation*}

It is possible to use the Volterra composition to express the solution of the following differential equation for every initial time $s \in I$.
\begin{equation}\label{eq:ode}
    \frac{d}{dt}y_s(t) =\tilde{f}(t) y_s(t), \quad y_s(s)=1 , \quad t \in [s, T ] \subset \mathbb{R};
\end{equation}
see, e.g., \cite{Giscard2015}.
 In fact, using Picard iterations, we get
\begin{align*}
\frac{d}{dt}y_s(t) &=\tilde{f}(t) y_s(t), \quad y_s(s)=1 \\
   & \downarrow \, \textrm{ integration } \\
y_s(t) &= 1 + \int_s^t \tilde{f}(\tau) y_s(\tau) \textrm{d}\tau \\
   & \downarrow \, \textrm{ integration } \\
y_s(t) &= 1 + \int_s^t \tilde{f}(\tau) \left(1 + \int_s^\tau \tilde{f}(\rho)y_s(\rho) \textrm{d}\rho \right) \textrm{d}\tau \\
 &= 1 + \int_s^t \tilde{f}(\tau) + \int_s^\tau \tilde{f}(\tau)\tilde{f}(\rho)y_s(\rho) \, \textrm{d}\rho \, \textrm{d}\tau  \\
    & \downarrow \, \textrm{ \dots } \\
y_s(t) &= 1 + \int_s^t \tilde{f}(\tau) \textrm{d}\tau  + \int_s^t \tilde{f}^{\star_v 2}(\tau) \textrm{d}\tau  + \dots ,
\end{align*}
from which we obtain the expression
\begin{equation}\label{eq:volterra:expression}
    y_s(t) = 1 + \int_s^t \sum_{k=1}^\infty \tilde{f}^{\star_v k}(\tau)\; \textrm{d}\tau.
\end{equation}

The Volterra composition is not a product and lacks essential features, for instance, the identity. For this reason, the Volterra composition has been extended, obtaining the so-called $\star$-product \cite{ProceedingsPaper2020}  that we briefly introduce in the following.
Consider the class $\D$ of all the distributions $d$ that can be written as
\begin{equation*}
d(t,s)=\widetilde{d}(t,s)\Theta(t-s) + \sum_{i=0}^N \widetilde{d}_i(t,s)\delta^{(i)}(t-s),  
\end{equation*}
where $N$ is a finite integer, $\tilde{d}, \tilde{d}_i$ are smooth bivariate functions over $I \times I$, $\Theta(\cdot)$ stands for the Heaviside theta function
\begin{equation*}
    \Theta(t-s) = \begin{cases}
                        1, \quad t \geq s \\
                        0, \quad t < s 
                \end{cases},
\end{equation*}
and $\delta^{(i)}(\cdot)$ is the $i$th derivative of the Dirac delta distribution $\delta(\cdot)=\delta^{(0)}(\cdot)$. 
We can endow the class $\D$ with a non-commutative algebraic structure by defining the $\star$-product as
\begin{equation}\label{eq:def:star}
  \big(f_2 \star f_1\big)(t,s) := \int_I f_2(t,\tau) f_1(\tau, s) \, \text{d}\tau, \quad f_1, f_2 \in \D.
\end{equation}
The $\star$-product is associative over $\D$, $\D$ is closed under $\star$-multiplication, and the identity element with respect to the $\star$-product is the Dirac delta distribution, $1_\star:=\delta(t-s)$, see, e.g., \cite{ProceedingsPaper2020}.

Consider the subclass $\Sm(I) \subset \D$ comprising those distributions of form
\begin{equation*}
f(t,s)=\widetilde{f}(t,s)\Theta(t-s).
\end{equation*}
Then, the $\star$-product between $f_1,f_2 \in \Sm(I)$ reduces to a Volterra composition
\begin{align*}
  \big(f_2 \star f_1\big)(t,s) &= \int_{I} \widetilde{f}_2(t,\tau) \widetilde{f}_1(\tau, s)\Theta(t-\tau)\Theta(\tau-s) \, \text{d}\tau,\\ &=\Theta(t-s)\int_s^{t} \widetilde{f}_2(t,\tau) \widetilde{f}_1(\tau, s) \, \text{d}\tau = \Theta(t-s)( \tilde{f}_2 \star_v \tilde{f}_1).
\end{align*}
As a consequence, using \eqref{eq:volterra:expression}, we can express the solution of \eqref{eq:ode} for every $s \in I$ as
\begin{equation}\label{eq:star:sol}
   y_s(t) = u(t,s) = \Theta(t-s) \star R_\star(f),     
\end{equation}
where $f(t,s) = \tilde{f}(t)\Theta(t-s)$ and $R_\star(f)$ is the $\star$-resolvent of $f$, i.e.,
\begin{equation*}
    R_\star(f) = \delta(t-s) + \sum_{k=1}^\infty f(t,s)^{\star k},
\end{equation*}
with $f(t,s)^{\star k} = \Theta(t-s) \tilde{f}(t)^{\star_v k}$. Note that the series $\sum_{k=1}^\infty f(t,s)^{\star k}$ converges for every $f \in \Sm(I)$.
The $\star$-product easily extends to matrices composed of elements from $\D$ by extending the scalar multiplication appearing in the integrand in \eqref{eq:def:star} to the usual matrix-matrix multiplication; see \cite{GiscardPozza2022} for more details.

While expression \eqref{eq:star:sol} is compact, the $\star$-resolvent definition hides an infinite series of nested integrals. Therefore, at first sight, it does not seem like a convenient expression. In the next section, we effectively solve this problem by showing that it is possible to approximate the $\star$-product by the usual matrix-matrix product. Consequently, for a fixed $s$, expression \eqref{eq:star:sol} can be approximated relatively cheaply by solving a linear system.

\section{Discretization of the $\star$-product}\label{sec:starsol}
		In this section, we describe an effective strategy for approximating the $\star$-product.
		Consider a sequence of orthonormal functions $\{p_k\}_k$ over the bounded interval $I=[0,T]$, i.e., 
		\begin{align*}
			\int_{I} p_k(\tau) p_\ell(\tau) d\tau = \begin{cases}
				0,\quad \text{if }k\neq \ell\\
				1,\quad \text{if } k=\ell				
			\end{cases},
		\end{align*}
		so that  $\{p_k\}_k$ is a basis for the space of smooth functions over $I$.
		Note that the functions $p_k$ are not in $\D$; hence we cannot (formally) $\star$-multiply them.
		Consider a distribution $f \in \Sm(I)$. The function $f(t,s)= \tilde{f}(t,s) \Theta(t-s)$ is piecewise smooth, therefore, we can choose the basis $\{p_k\}_k$ so that
        \begin{equation}\label{eq:f:exp}
			f(t,s) = \sum_{k=0}^\infty \sum_{\ell=0}^\infty f_{k,\ell} \, p_k(t) p_\ell(s), \quad t \neq s, \quad t,s \in I,
		\end{equation}
		with coefficients
		\begin{equation*}
		    f_{k,\ell} = \int_I \int_I f(\tau,\rho) p_k(\tau) p_\ell(\rho) \; \textrm{d} \rho \; \textrm{d} \tau.
		\end{equation*}
		For instance, the basis $\{p_k\}_k$ can be set as the sequence of shifted Legendre polynomials (e.g., \cite[p.~55]{LebSil72}).
		Defining the \emph{coefficient matrix} $F_M$ and the vector $\phi_M(t)$ as
		\begin{equation}\label{eq:coeff:mtx}
		          F_M := \begin{bmatrix}
				f_{0,0} & f_{0,1} & \dots & f_{0,M-1}\\
				f_{1,0} & f_{1,1} & \dots & f_{1,M-1}\\
				\vdots & \vdots &  & \vdots\\
				f_{M-1,0} & f_{M-1,1} & \dots & f_{M-1,M-1}
		\end{bmatrix}, 
		\quad  \phi_M(t) :=
		 \begin{bmatrix}
			p_0(s)\\
			p_1(s)\\
			\vdots\\
			p_{M-1}(s)
		\end{bmatrix},
		\end{equation}
		the truncated expansion series can be written in the matrix form:
	\begin{align*}
	f_M(t,s) := \sum_{k=0}^{M-1} \sum_{\ell=0}^{M-1} f_{k,\ell} \, p_k(t) p_\ell(s) = \phi_M(t)^T F_M \, \phi_M(s).
		\end{align*}
  
  Consider $f,g,h \in \Sm(I)$ so that $h = f \star g$, and the related coefficient matrices \eqref{eq:coeff:mtx}, respectively, $F_M,G_M,H_M$. By replacing $f$ and $g$ with their expansion \eqref{eq:f:exp}, it is not difficult to show that the expansion coefficients for $h$ are given by
  \begin{equation}\label{eq:hcoeff}
      h_{k,\ell} = \sum_{j=0}^{\infty} f_{k,j} \, g_{j,\ell}.
  \end{equation}
  As a consequence, we can approximate $H_M$ by the expression
  \begin{equation}\label{eq:mtx:prod}
       H_M \approx \hat{H}_m := F_M G_M,
  \end{equation}
  i.e., the $\star$-product can be approximated by the usual matrix-matrix multiplication of the related coefficient matrices.
  
  The approximation \eqref{eq:mtx:prod} is affected by a truncation error. Therefore, fixing $k$ and $\ell$, if the magnitude of the product $f_{k,j}, g_{j,\ell}$ in \eqref{eq:hcoeff} does not decay quickly enough for $j \rightarrow \infty$, then the truncation error $(H_M)_{k,\ell} - (\hat{H}_M)_{k,\ell}$ can be too large for practical purposes. Luckily, since $f \in \Sm$, numerical considerations illustrate that $F_M$ and $G_M$ are numerically banded for a certain choice of $\{ p_k \}_k$; for instance, see  Section \ref{sec:matrices} where we choose the shifted Legendre polynomials. Therefore, $M$ does not need to be too large to reach a small truncation error in the approximation \eqref{eq:mtx:prod}, excluding the last rows of the matrix $\hat{H}_M$ where the truncation error can still be significant. Further details and explanations on this matter are being developed and will be presented in future work. For the moment, in Section \ref{sec:matrices}, we provide numerical evidence of these claims.
  
  To conclude the presentation, we must discuss the convergence behavior of expansion \eqref{eq:f:exp}. Indeed, since $f$ is discontinuos for $t=s$, the expansion may not converge quickly (or may not converge) to $f(t,s)$ for every $t,s \in I$; see, e.g., \cite{LebSil72,TreApp13} for the polynomial case. Nevertheless, fixing $s = 0$, the univariate function $f(t,0) = \tilde{f}(t,0)\Theta(t-0) = \tilde{f}(t,0)$ is smooth over $I=[0,T]$. Therefore
  \begin{equation*}
      	f(t,0) = \sum_{k=0}^\infty a_k p_k(t) = \sum_{k=0}^\infty  p_k(t) \sum_{\ell=0}^\infty f_{k,\ell} \, p_\ell(0),
  \end{equation*}
  with $a_k = \sum_{\ell=0}^\infty (f_{k,\ell}\, p_\ell(0))$. As a consequence, we can approximate the function $f(t,0)$ by the expression
  \begin{equation*}
      f(t,0) \approx \phi_M(t)^T F_M \, \phi_M(0),
  \end{equation*}
  and expect to reach a small enough accuracy for a (relatively) small $M$.
 Section~\ref{sec:numexp} illustrates with several numerical examples that it is possible to achieve machine precision accuracy for a small value of $M$. 
  
   Consider the function $u(t,s)$ in \eqref{eq:star:sol}. Using the previous construction, the related coefficient matrix $U_M$, i.e., such that $u(t,s) \approx \phi(t)_M^T\, U_M\, \phi_M(s)$, can be approximated by
   \begin{equation*}
       U_M \approx T_M (I_M-F_M)^{-1},
   \end{equation*}
   where $T_M$ is the coefficient matrix of $\Theta(t-s)$, and $F_M$ is the coefficient matrix of $\tilde{f}(t)\Theta(t-s)$, with $\tilde{f}$ from \eqref{eq:ode}.
   Since $u(t,s) \in \Sm$, for $s=0$ we can approximate the solution of \eqref{eq:ode} by the formula:
   \begin{equation*}
       y_0(t) \approx \phi_M(t)^T U_M \, \phi_M(0) \approx \phi_M(t)^T T_M (I_M-F_M)^{-1}  \phi_M(0).
   \end{equation*}
  Then, the vector $u_M = T_M  x$ contains the approximated expansion coefficients of $y_0(t)$, i.e.,
  \begin{equation*}
      y_0(t) \approx \phi_M(t)^T u_M, \; \text{ for every } t \in I,
  \end{equation*}
  where $x$ is the solution of the linear system
  \begin{equation}\label{eq:linsys}
      (I_M-F_M) x =  \phi_M(0).
  \end{equation}

\section{Properties of the coefficient matrix}\label{sec:matrices}
In this section, we illustrate several properties of the coefficient matrices \eqref{eq:coeff:mtx} through numerical examples. We set $I = [0,1]$, and, as the sequence of orthonormal functions, we choose the sequence of orthonormal shifted Legendre polynomials, i.e., the sequence of polynomials $\{p_k\}_k$ such that
		\begin{align*}
			\int_{0}^1 p_k(\tau) p_\ell(\tau) d\tau = \begin{cases}
				1,\quad \text{if } k=\ell	\\
				0,\quad \text{if }k\neq \ell
			\end{cases},
		\end{align*}
with $k, \ell$ the degree of the polynomial.
In the following, we consider the functions $f_k(t,s) = \tilde{f}_k(t)\Theta(t-s)$ from Table \ref{table:band} and the related $M \times M$ coefficient matrix $F^{(k)}_M$ defined in \eqref{eq:coeff:mtx}. The numerical experiments were performed using MatLab R2022a.

\begin{table}[t]
\begin{center}
\begin{tabular}
[c]{lccccc}
\hline
\rule{0pt}{13pt}
Functions             &  $\tilde{f}_1 = 1$    & $\tilde{f}_2 = t$   & $\tilde{f}_3 = t^3$   & $\tilde{f}_4 = \cos(t)$   &  $\tilde{f}_5 = \log(t+1) $ \\
\hline
\multicolumn{6}{c}{ $M = 25$ }\\
\hline
 Num. band.            & $1$ & $2$  & $4$  & $13$ & $20$  \\
 Spectral radius           & $0.0592$ & $0.0357$  & $0.0238$  & $0.0480$ & $0.0271$  \\
 $\sigma_{\min}$           & $2.42e-3$ & $3.50e-5$  & $9.68e-9$  & $1.74e-3$ & $3.42e-5$  \\
 $\sigma_{\max}$           & $1.2732$ & $0.9447$  & $0.6864$  & $0.9694$ & $0.6938$  \\
\hline
\multicolumn{6}{c}{ $M = 100$ }\\
 \hline
 Num. band.            & $1$ & $2$  & $4$  & $13$ & $20$  \\
 Spectral radius           & $0.0556$ & $0.0296$  & $0.0155$  & $0.0444$ & $0.0223$  \\
 $\sigma_{\min}$           & $1.56e-4$ & $1.57e-7$  & $2.33e-13$  & $1.10e-4$ & $1.56e-7$  \\
 $\sigma_{\max}$           & $1.2732$ & $0.9447$  & $0.6864$  & $0.9694$ & $0.6938$  \\
\hline
\multicolumn{6}{c}{ $M = 500$ }\\
 \hline
 Num. band.            & $1$ & $2$  & $4$  & $13$ & $20$  \\
 Spectral radius           & $0.0554$ & $0.0297$  & $0.0146$  & $0.0458$ & $0.0215$  \\
 $\sigma_{\min}$           & $6.27e-6$ & $2.59e-10$  & $6.58e-19$  & $2.59e-10$ & $3.42e-5$  \\
 $\sigma_{\max}$           & $1.2732$ & $0.9447$  & $0.6864$  & $0.9694$ & $0.6938$  \\
 \hline
\end{tabular}
\end{center}
\caption{\label{table:band} Properties of coefficient matrices $F^{(k)}_M$.}
\end{table}

Table \ref{table:band} reports the numerical bandwidth of each coefficient matrix $F_{M}^{(k)}$ for $M=25,100, 500$. With numerical bandwidth, we mean the bandwidth of the matrix once all its elements with a magnitude smaller than the machine precision have been rounded to zero. 
First, we observe that the numerical bandwidth is the same for every value of $M$. Moreover, we note that for the polynomial functions $\tilde{f}_1, \tilde{f}_2, \tilde{f}_3$, the corresponding bandwidth is equal to the degree of the polynomial plus one. Finally, the functions $f_4(t,s)=\cos(t)\Theta(t-s), f_5(t,s)= \log(t+1)\Theta(t-s)$ are also numerically banded.

\begin{figure}[p]
\begin{center}
\includegraphics[width=6.5cm]{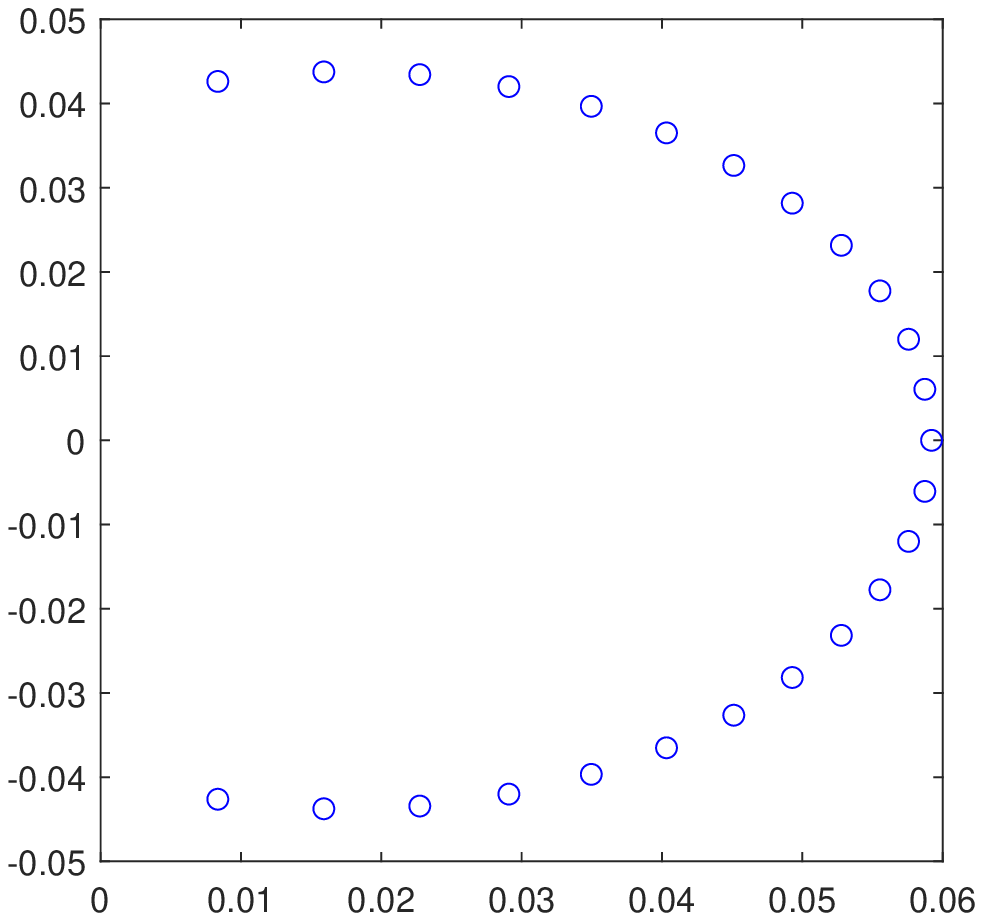}
\includegraphics[width=6.5cm]{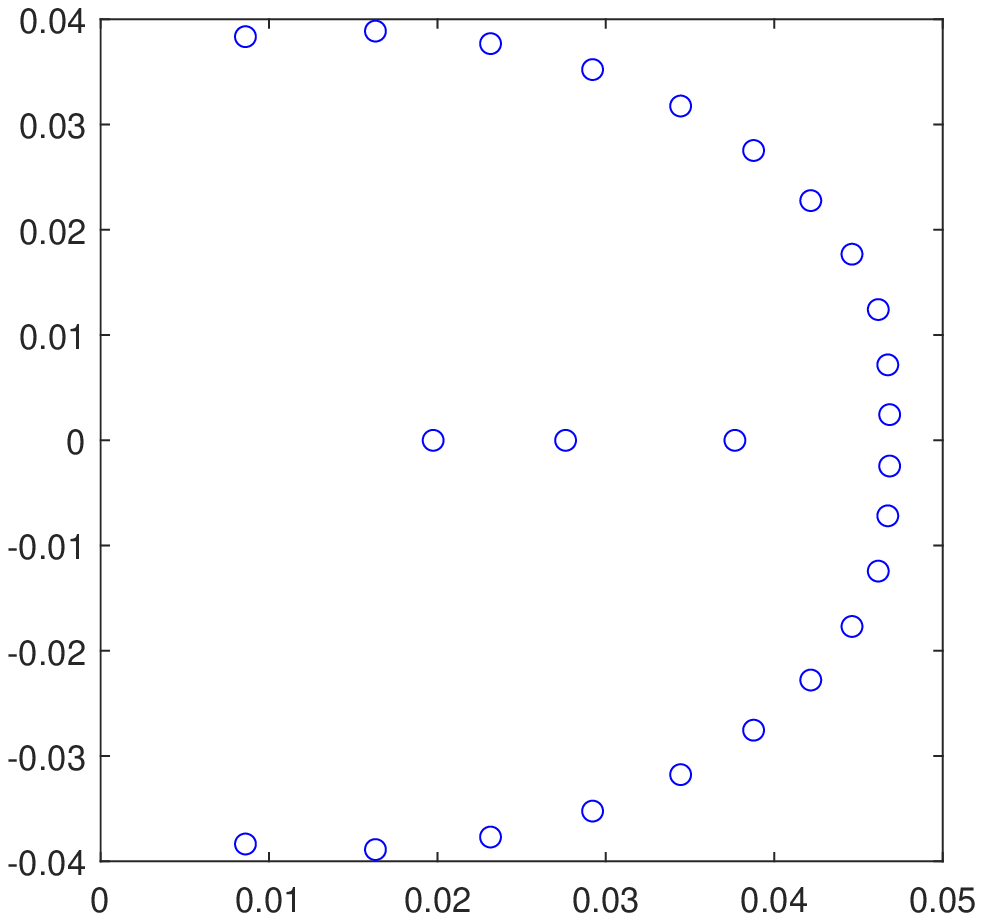} \\
{\small $M=25$ }\\
\includegraphics[width=6.5cm]{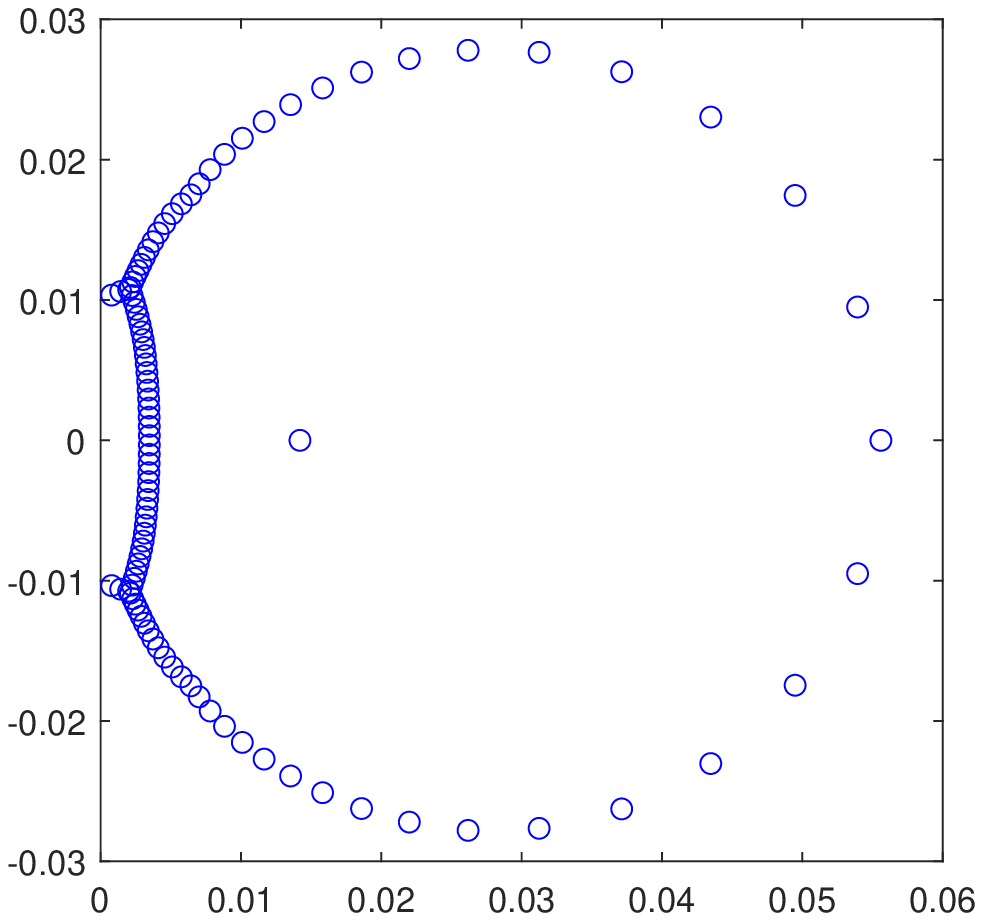}
\includegraphics[width=6.5cm]{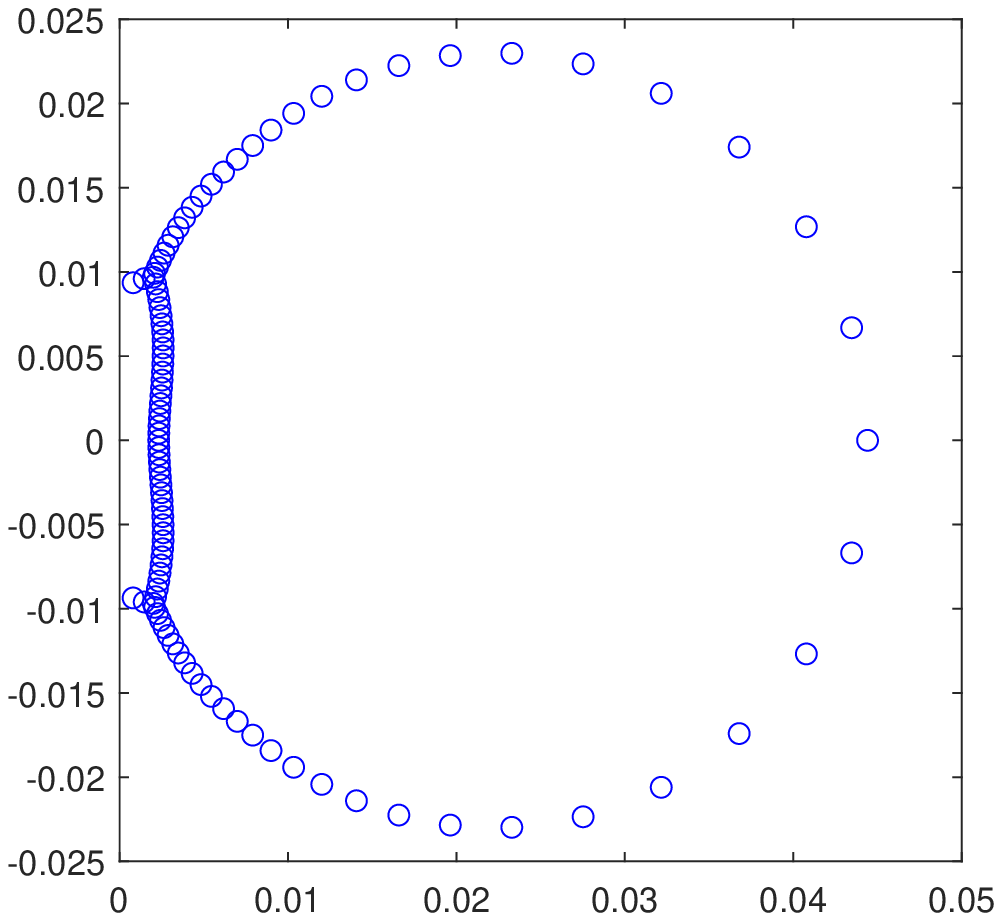} \\
{\small $M=100$ }\\
\includegraphics[width=6.5cm]{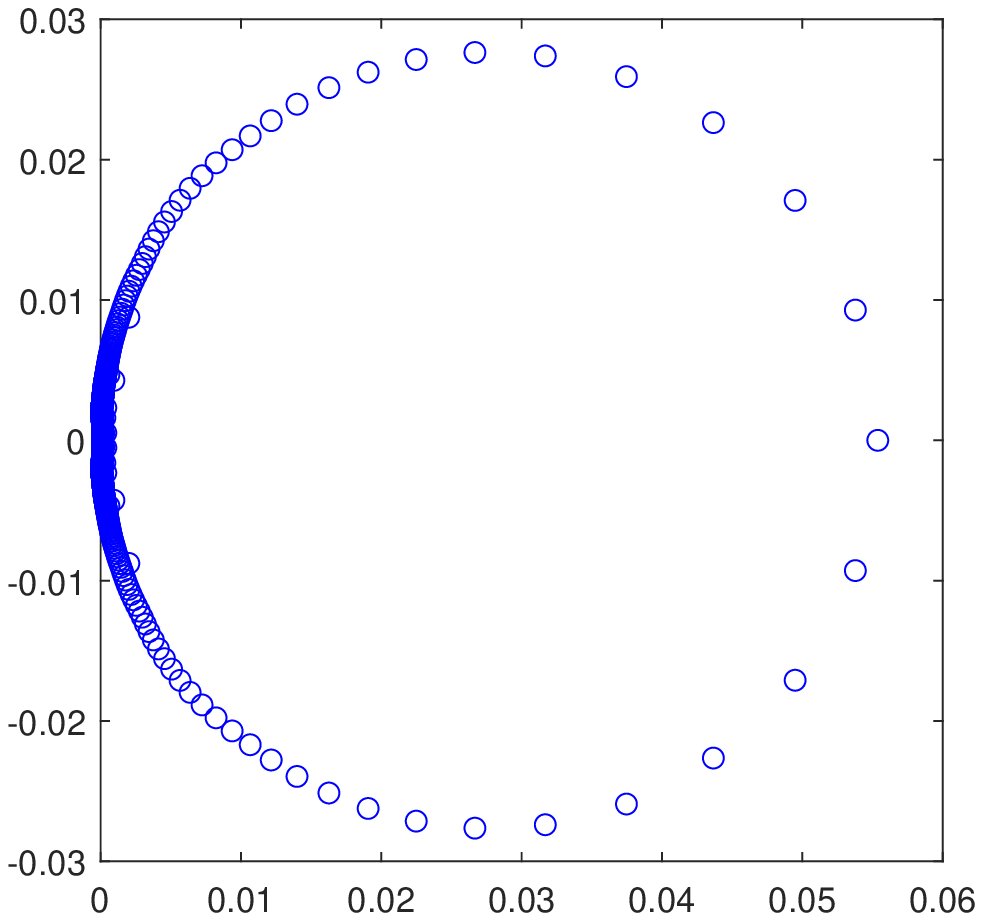}
\includegraphics[width=6.5cm]{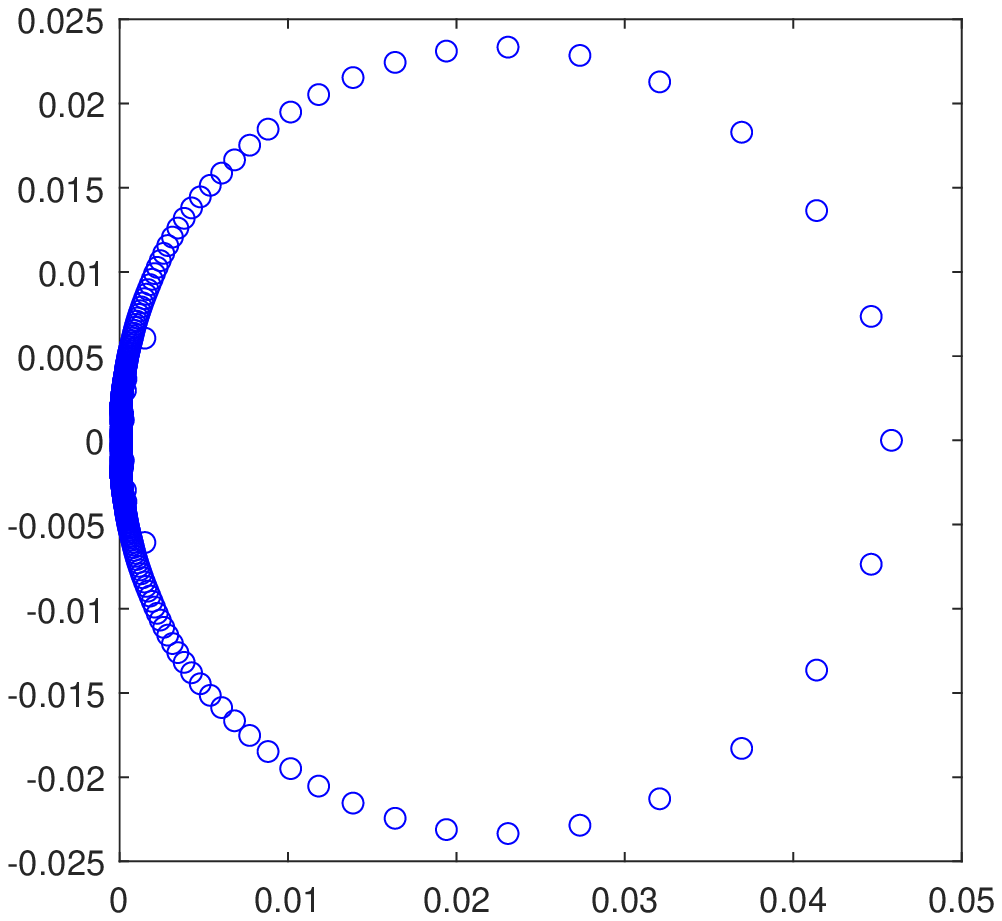} \\
{\small $M=500$ }\\
\caption{\label{fig:spectrum} Spectrum of the coefficient matrices $F_M^{(1)}$ (left), $F_M^{(4)}$ (right) defined in \eqref{eq:coeff:mtx}, for $M= 25, 100, 500$.}
\end{center}
\end{figure}

Table \ref{table:band} also reports the spectral radius and the minimal and maximal singular values (respectively $\sigma_{\min}$, $\sigma_{\max}$) of each $F_M^{(k)}$. While both the spectral radius and $\sigma_{\max}$ do not vary significantly for $M=25, 100, 500$, $\sigma_{\min}$ becomes smaller as $M$ increases. As the linear system \eqref{eq:linsys} involves the shifted matrix $I_M - F_M^{(k)}$, it is important to note that all the computed spectral radii are smaller than $1$. 

Finally, Figure \ref{fig:spectrum} presents the spectra of the matrices $F_M^{(1)}$ and $F_M^{(4)}$ for $M=25, 100, 500$. For both the functions, as $M$ increases, the spectrum tends to distribute in a circle on the right-half of the complex plane, closer and closer to the origin. We do not report the spectrum plots of the other matrices considered above since they display analogous behavior.

\section{Numerical experiments}\label{sec:numexp}
In this section, we test the numerical method explained in Section \ref{sec:starsol} on the ODE
\begin{equation}\label{eq:ode:experiments}
    \frac{d}{d t} y(t) = \tilde{f}_k(t) y(t), \quad y(0) = 1, \quad t \in I = [0,1],
\end{equation}
for each function $\tilde{f}_k$ from Table \ref{table:band}.
More precisely, the method works as follows:
\begin{enumerate}
    \item We discretize $f_k(t,s)= \tilde{f}_k(t) \Theta(t-s)$ as described in Section \ref{sec:starsol}, obtaining the matrix $F_M^{(k)}$. We use as an orthonormal basis the shifted orthonormal Legendre polynomials in Section \ref{sec:matrices}.
    \item Let  $b$ be the numerical bandwidth of $F_M^{(k)}$; we define the matrix $\hat{F}_M^{(k)}$ by setting the last $b$ rows of the matrix $F_M^{(k)}$ to zero. This has proven helpful in reducing the accumulation of truncation errors in the last rows of the solution.
    \item We solve the (banded) linear system
\begin{equation*}
    \left (I_M - \hat{F}_M^{(k)} \right) x = \phi_M(0),
\end{equation*}
using the Matlab backslash \texttt{\symbol{92}} operation.
    \item The solution of \eqref{eq:ode:experiments} is given by
    \begin{equation}\label{eq:starsol:experiments}
    y(t) \approx \hat{y}_M(t) := \phi_M(t)^T u_M, \quad \textrm{ with } u_M = T_m \, x.
\end{equation}
 \end{enumerate}
 
In Table \ref{table:exp}, we report the maximal relative error of approximation \eqref{eq:starsol:experiments} over $I$ for $M=25, 100$. The relative errors were computed on an equispaced mesh of $100$ points over $[0,1]$. 
As a reference value for the solution, we considered the function $\exp(\int_0^t \tilde{f}_k(\tau) \, \text{d}\tau)$.
We compare our results with the maximal relative errors obtained using the Matlab methods \texttt{ode45} and \texttt{ode89} with relative and absolute tolerances set equal to $ \texttt{eps} = 2.2204e-16$.
 For $M=100$, Table \ref{table:exp} shows that approximation \eqref{eq:starsol:experiments} is always better than the others.  On the other hand, for $M=25$, we obtain worse results for $k=3,4,5$, showing that it is possible to calibrate the accuracy of the solution by the matrix size.

With these experiments, we do not want to claim anything about the performance of our method compared to well-established explicit methods such as  \texttt{ode45} and \texttt{ode89}. The examples considered here are certainly not enough for drawing any conclusion. The table aims to show that approximation \eqref{eq:starsol:experiments} can compete in accuracy with well-established approaches, a promising result for our future work.

\begin{table}[ht]
\begin{center}
\begin{tabular}
[c]{lccccc}
\hline
\rule{0pt}{13pt}
Functions             &  $\tilde{f}_1 = 1$    & $\tilde{f}_2 = t$   & $\tilde{f}_3 = t^3$   & $\tilde{f}_4 = \cos(t)$   &  $\tilde{f}_5 = \log(t+1) $ \\
\hline
 $\hat{y}_{25}(t)$  & $1.20e-15$ & $1.11e-15$ & $3.36e-14$  & $1.37e-09$ & $4.04e-04
$  \\
 $\hat{y}_{100}(t)$ & $1.20e-15$ & $1.11e-15$  & $8.88e-16$  & $1.22e-15$ & $9.77e-16$  \\
 \texttt{ode45}            & $9.76e-15$ & $4.61e-13$  & $1.53e-12$  & $1.13e-13$ & $9.72e-13$  \\
 \texttt{ode89}            & $1.17e-13$ & $6.69e-14$  & $3.63e-14$  & $1.13e-13$ & $9.92e-14$  \\
\hline
\end{tabular}
\end{center}
\caption{\label{table:exp} Maximal relative error over $I=[0,1]$ of \texttt{ode45}, \texttt{ode89} methods and of the approximation $\hat{y}_M(t)$ in \eqref{eq:starsol:experiments} for the solution of the ODE \eqref{eq:ode:experiments}.}
\end{table}

\section{Conclusion and future work}\label{sec:conclusion}
In this work, we have explained how to express the solution of a scalar linear ODE using the so-called $\star$-product. Moreover, we have shown how to derive a numerical method from this expression and successfully tested it on several examples. The numerical method requires solving a linear system whose properties have also been numerically investigated.
Concerning the numerical efficiency of the introduced method, other possible approaches in the solution of the linear system may be used -- for instance, Krylov subspace methods. Furthermore, since the solution depends continuously on the initial time parameter $s$, we are also investigating the use of acceleration methods such as the one in \cite{BuoKarPoz15}. In addition, we are currently developing an efficient method for computing the coefficient matrix $F_m$.

 Given a smooth matrix-valued function $\tilde{A}(t) \in \mathbb{C}^{N \times N}$, the solution of the system
\begin{equation*}
    \frac{d}{d t} Y_s(t) = \tilde{A}(t) Y_s(t), \quad y(s) = I_N, \quad t\geq s, \quad t,s \in I,
\end{equation*}
can also be expressed as
\begin{equation*}
    Y_s(t) = U(t,s) = \Theta(t-s) \star R_{\star}(\tilde{A}(t)\Theta(t-s)), \quad t,s \in I,
\end{equation*}
following the results in \cite{Giscard2015}.
Therefore, the scalar method we have described can be generalized to the more challenging problem offered by systems of non-autonomous linear ODEs. The results discussed in this work are thus promising for developing new efficient methods for computing $Y_s(t)$.

\section*{Acknowledgements}
This  work  was  supported  by  Charles  University  Research programs UNCE/SCI/023 and PRIMUS/21/SCI/009 and by the Magica project ANR-20-CE29-0007 funded by the French National Research Agency.


\end{document}